\documentclass[a4paper,10pt,english]{smfart}
\usepackage[english]{babel}
\usepackage[OT2,T1]{fontenc}
\usepackage{lmodern}
\usepackage{graphicx}
\usepackage{tabularx}
\usepackage{amsmath,amsthm,amssymb,amsfonts}
\usepackage[a4paper,left=4cm,right=4cm,top=4cm,bottom=4cm]{geometry}
\numberwithin{equation}{section}
\title[Distribution of squarefree integers in arithmetic progressions]{On the distribution of squarefree integers in arithmetic progressions}
\author{Pierre Le Boudec}
\subjclass{$11$N$37$, $11$N$69$}
\address{EPFL SB MATHGEOM TAN \\ MA C$3$ $604$ (B\^{a}timent MA) \\ Station $8$ \\ \text{CH-$1015$} Lausanne \\ Switzerland}
\email{pierre.leboudec@epfl.ch}

\begin{document}

\makeatletter
\def\imod#1{\allowbreak\mkern10mu({\operator@font mod}\,\,#1)}
\makeatother

\newtheorem{lemma}{Lemma}
\newtheorem{theorem}{Theorem}
\newtheorem{corollary}{Corollary}
\newtheorem{proposition}{Proposition}
\newtheorem{conjecture}{Conjecture}
\newtheorem{conj}{Conjecture}
\renewcommand{\theconj}{\Alph{conj}}

\newcommand{\vol}{\operatorname{vol}}
\newcommand{\D}{\mathrm{d}}
\newcommand{\rank}{\operatorname{rank}}
\newcommand{\Pic}{\operatorname{Pic}}
\newcommand{\Gal}{\operatorname{Gal}}
\newcommand{\meas}{\operatorname{meas}}
\newcommand{\Spec}{\operatorname{Spec}}
\newcommand{\eff}{\operatorname{eff}}
\newcommand{\rad}{\operatorname{rad}}
\newcommand{\sq}{\operatorname{sq}}
\newcommand{\tors}{\operatorname{tors}}
\newcommand{\Cl}{\operatorname{Cl}}

\begin{abstract}
We investigate the error term of the asymptotic formula for the number of squarefree integers up to some bound, and lying in some arithmetic progression $a \imod{q}$. In particular, we prove an upper bound for its variance as $a$ varies over
$(\mathbb{Z}/q\mathbb{Z})^{\times}$ which considerably improves upon earlier work of Blomer.
\end{abstract}

\maketitle

\tableofcontents

\section{Introduction}

\subsection{Distribution in arithmetic progressions: analogy between primes and squarefree integers}

Let $x \geq 1$ and $q \in \mathbb{Z}_{\geq 1}$ be such that $q \leq x$, and let also
$a \in (\mathbb{Z} / q \mathbb{Z})^{\times}$. We introduce
\begin{equation*}
S(x;q,a) = \sum_{\substack{n \leq x \\ n = a \imod{q}}} |\mu(n)|,
\end{equation*}
and we define the quantity $E(x;q,a)$ by the equality
\begin{equation*}
S(x;q,a) = c_q \frac{x}{q} + E(x;q,a),
\end{equation*}
where
\begin{equation*}
c_q = \prod_{p \nmid q} \left( 1 - \frac1{p^2} \right).
\end{equation*}
Hooley \cite[Theorem $3$]{MR0371799} proved that for any fixed $\varepsilon > 0$, we have
\begin{equation}
\label{Hooley}
E(x;q,a) \ll \left( \frac{x}{q} \right)^{1/2} + q^{1/2 + \varepsilon},
\end{equation}
which is non-trivial provided that $q \leq x^{2/3 - \varepsilon}$, and thus gives an asymptotic formula for $S(x;q,a)$ in this range. It should be stated that it is expected to be difficult to improve on this result.

It is thus natural to study this problem on average. More specifically, the purpose of this article is to investigate the following variance
\begin{equation*}
V(x;q) = \sideset{}{^\ast}\sum_{a \imod{q}} E(x;q,a)^2,
\end{equation*}
where the symbol $\ast$ indicates that the summation is restricted to residue classes
$a \in (\mathbb{Z} / q \mathbb{Z})^{\times}$. Blomer studied this quantity and he proved
\cite[Theorem $1$.$3$]{MR2444061} that
\begin{equation}
\label{Blomer}
V(x;q) \ll x^{\varepsilon} \left( x + \min \left\{ \frac{x^{5/3}}{q},  q^2 \right\} \right).
\end{equation}
Unfortunately, this upper bound is not sharp (even without the factor $x^{\varepsilon}$) unless the size of $q$ is about
$x$. Croft established a Barban-Davenport-Halberstam type Theorem \cite[Theorem $2$]{MR0379402}, which states that $V(x;q)$ has size about $x^{1/2} q^{1/2}$ on average over $q \leq Q$, and for $Q \geq x^{2/3}$. Note that Croft was summing over all residue classes $a \in \mathbb{Z} / q \mathbb{Z}$ but this can certainly be ignored in our discussion.

Croft noticed that his result is in accordance with a remark of Montgomery stating that on probabilistic grounds, one would expect
\begin{equation*}
E(x;q,a) \ll \left( \frac{x}{q} \right)^{1/4 + \varepsilon}.
\end{equation*}
However, it should be pointed out that the range for which this upper bound should hold was unspecified. In particular, it is not clear why this should be true if $q$ is very close to $x$.

It is very instructive to compare the problem at hand with the analog problem about primes in arithmetic progressions. Let us define the quantity $\mathcal{E}(x;q,a)$ by the equality
\begin{equation*}
\sum_{\substack{n \leq x \\ n = a \imod{q}}} \Lambda(n) = \frac{x}{\varphi(q)} + \mathcal{E}(x;q,a),
\end{equation*}
and let us also introduce the corresponding variance
\begin{equation*}
\mathcal{V}(x;q) = \sideset{}{^\ast}\sum_{a \imod{q}} \mathcal{E}(x;q,a)^2.
\end{equation*}
After proving that an expectation of Montgomery \cite{MR0427249} was false, Friedlander and Granville
\cite[Conjecture $1$]{MR986796} conjectured the following.

\begin{conj}
\label{Friedlander-Granville}
Let $\varepsilon > 0$ be fixed. For $q \leq x$, we have
\begin{equation*}
\mathcal{E}(x;q,a) \ll \left( \frac{x}{q} \right)^{1/2} x^{\varepsilon}.
\end{equation*}
\end{conj}

Note that if $q \geq x^{\delta}$ for some $\delta > 0$, this conjecture is stronger than the Generalized Riemann Hypothesis (GRH), which only gives
\begin{equation}
\label{Consequence GRH Primes}
\mathcal{E}(x;q,a) \ll x^{1/2} (\log x)^2.
\end{equation}
However, Tur\'{a}n \cite{zbMATH03026443} proved that GRH implies that for $q \leq x$, we have
\begin{equation}
\label{Turan}
\mathcal{V}(x;q) \ll x (\log x)^4,
\end{equation}
which agrees with Conjecture \ref{Friedlander-Granville}.

We now make the following conjecture, which is the analog of Conjecture \ref{Friedlander-Granville} for our present problem.

\begin{conjecture}
\label{Conjecture}
Let $\varepsilon > 0$ be fixed. For $q \leq x$, we have
\begin{equation*}
E(x;q,a) \ll \left( \frac{x}{q} \right)^{1/4} x^{\varepsilon}.
\end{equation*}
\end{conjecture}

As in the case of primes in arithmetic progressions, this conjecture can be analyzed by looking at the poles and residues of the complex functions $L(\chi,s)/L(\chi^2,2s)$ where $\chi$ runs over the set of Dirichlet characters modulo $q$. It is worth pointing out that here, similarly to \eqref{Consequence GRH Primes}, it should not be expected that GRH implies anything stronger (maybe up to the factor $x^{\varepsilon}$) than
\begin{equation}
\label{Consequence GRH}
E(x;q,a) \ll x^{1/4 + \varepsilon} q^{1/4}.
\end{equation}
Note that Hooley's upper bound \eqref{Hooley} is stronger than \eqref{Consequence GRH} for
$q \geq x^{1/3 - \varepsilon}$, but going below this range seems to be challenging. In particular, the upper bound \eqref{Consequence GRH} for $q=1$ implies the Riemann Hypothesis.

Surprisingly, Moreira Nunes \cite[Corollary $1$.$3$]{MoreiraNunes} recently established, for $q$ large enough in terms of $x$, a more precise and unconditional analog of Tur\'{a}n's result \eqref{Turan}. More specifically, he proved that for
$x^{31/41 + \varepsilon} \leq q \leq x^{1 - \varepsilon}$, we have
\begin{equation}
\label{Moreira Nunes}
V(x;q) \sim C_q x^{1/2} q^{1/2},
\end{equation}
where $C_q > 0$ is explicit and uniformly bounded. Note that this is in agreement with Conjecture \ref{Conjecture}.

It would be very interesting to have a heuristic giving the maximal range for which the asymptotic formula
\eqref{Moreira Nunes} should hold. It seems to the author that using the result of Heath-Brown \cite[Lemma $3$]{MR757475} (as in section \ref{Section Proposition}) instead of the square sieve in the work of Moreira Nunes, the error term in \cite[Theorem $1$.$1$]{MoreiraNunes} can be replaced by
$x^{\varepsilon} \left( x^{1/3} q^{2/3} + x^{5/3} q^{-4/3} + x^{4/3} q^{-2/3} \right)$. Therefore, it can be checked that the asymptotic formula \eqref{Moreira Nunes} actually holds for
$x^{5/7 + \varepsilon} \leq q \leq x^{1 - \varepsilon}$.

It is also natural to ask in which range one can prove the analog of the result of Tur\'{a}n \eqref{Turan}. As already stated above, the upper bound $V(x;1) \ll x^{1/2 + \varepsilon}$ implies the Riemann Hypothesis so this problem is extremely hard if
$q$ is small. Nevertheless, we can prove the following result.

\begin{theorem}
\label{Main Theorem}
Let $\varepsilon > 0$ be fixed. For $q \leq x$, we have
\begin{equation*}
V(x;q) \ll x^{\varepsilon} \left( x^{1/2} q^{1/2} + \frac{x}{q^{1/2}} \right).
\end{equation*}
In particular, for $x^{1/2} \leq q \leq x$, we have
\begin{equation*}
V(x;q) \ll x^{1/2 + \varepsilon} q^{1/2}.
\end{equation*}
\end{theorem}

Theorem \ref{Main Theorem} improves upon the upper bound of Blomer \eqref{Blomer} in the whole range $q \leq x$, and also on the upper bound of Moreira Nunes \cite[Theorem $1$.$1$]{MoreiraNunes} in the range where the asymptotic formula \eqref{Moreira Nunes} is not known to hold.

It should be noted that in the first upper bound in Theorem \ref{Main Theorem}, the only thing which could be improved without knowing a quasi-Riemann Hypothesis (apart from the factor $x^{\varepsilon}$) is the power of $q$ in the term $x q^{-1/2}$.

Theorem \ref{Main Theorem} will be deduced from the investigation of the quantity
\begin{equation*}
T(x;q) = \sideset{}{^\ast}\sum_{\substack{n_1, n_2 \leq x \\ n_1 = n_2 \imod{q}}} |\mu(n_1)| |\mu(n_2)|,
\end{equation*}
where, as for residue classes, the symbol $\ast$ indicates that the summation is restricted to integers which are coprime to $q$. We will prove that Theorem \ref{Main Theorem} is equivalent to the following asymptotic formula for $T(x;q)$.

\begin{proposition}
\label{Proposition}
Let $\varepsilon > 0$ be fixed. For $q \leq x$, we have
\begin{equation*}
T(x;q) - \frac1{\varphi(q)} \left( \sideset{}{^\ast}\sum_{n \leq x} |\mu(n)| \right)^2 \ll 
x^{\varepsilon} \left( x^{1/2} q^{1/2} + \frac{x}{q^{1/2}} \right).
\end{equation*}
\end{proposition}

It is worth noting that Proposition \ref{Proposition} is equivalent to the following upper bound for the variance of certain character sums. 

\begin{corollary}
\label{Corollary 0}
Let $\varepsilon > 0$ be fixed. For $q \leq x$, we have
\begin{equation*}
\frac1{\varphi(q)} \sum_{\chi \neq \chi_0} \left| \sum_{n \leq x} |\mu(n)| \chi(n) \right|^2 \ll x^{\varepsilon}
\left( x^{1/2} q^{1/2} + \frac{x}{q^{1/2}} \right),
\end{equation*}
where the first sum is over all non-trivial Dirichlet characters modulo $q$.
\end{corollary}

The proof of Proposition \ref{Proposition} uses results established by the author in \cite{MR3263143} to count solutions to certain linear congruences, and which draw upon geometry of numbers methods.

\subsection{Density results}

Several interesting density results immediately follow from Theorem \ref{Main Theorem}. First of all, we obtain that in the range
$x^{1/2} \leq q \leq x$, the set of $a  \in (\mathbb{Z} / q \mathbb{Z})^{\times}$ violating Conjecture \ref{Conjecture} has density zero.

\begin{corollary}
\label{Corollary 1}
Let $\varepsilon > 0$ be fixed. For $x^{1/2} \leq q \leq x$, we have
\begin{equation*}
\frac1{\varphi(q)} \# \left\{ a \in (\mathbb{Z} / q \mathbb{Z})^{\times}, |E(x;q,a)| > \left( \frac{x}{q} \right)^{1/4} x^{\varepsilon} \right\} \ll x^{- \varepsilon}.
\end{equation*}
\end{corollary}

Corollary \ref{Corollary 1} states that for $x^{1/2} \leq q \leq x^{1 - \varepsilon}$ and for almost every
$a \in (\mathbb{Z} / q \mathbb{Z})^{\times}$, we have an asymptotic formula for $S(x;q,a)$ in which the error term is about the fourth root of the main term, which is a strong improvement of Hooley's error term $q^{1/2 + \varepsilon}$ appearing in \eqref{Hooley}.

We have the following analog concerning the range $q \leq x^{1/2}$.

\begin{corollary}
\label{Corollary 2}
Let $\varepsilon > 0$ be fixed. For $q \leq x^{1/2}$, we have
\begin{equation*}
\frac1{\varphi(q)} \# \left\{ a \in (\mathbb{Z} / q \mathbb{Z})^{\times},
|E(x;q,a)| > \frac{x^{1/2 + \varepsilon}}{q^{3/4}} \right\} \ll x^{- \varepsilon}.
\end{equation*}
\end{corollary}

Corollary \ref{Corollary 2} states that in the range $q\leq x^{1/2} $ and for almost every
$a \in (\mathbb{Z} / q \mathbb{Z})^{\times}$, we have an asymptotic formula for $S(x;q,a)$ in which the error term is
$x^{1/2 + \varepsilon} q^{-3/4}$, which is a modest improvement of Hooley's error term $x^{1/2} q^{-1/2}$ appearing in \eqref{Hooley}.

Recall that Hooley's upper bound \eqref{Hooley} is only stronger than \eqref{Consequence GRH} in the range
$q \geq x^{1/3 - \varepsilon}$. Theorem \ref{Main Theorem} implies that for $x^{1/4} \leq q \leq x^{1/3}$, the set of
$a \in (\mathbb{Z} / q \mathbb{Z})^{\times}$ violating the upper bound \eqref{Consequence GRH} has density zero.

\begin{corollary}
\label{Corollary 3}
Let $\varepsilon > 0$ be fixed. For $x^{1/4} \leq q \leq x^{1/3}$, we have
\begin{equation*}
\frac1{\varphi(q)} \# \left\{ a \in (\mathbb{Z} / q \mathbb{Z})^{\times}, |E(x;q,a)| > x^{1/4 + \varepsilon} q^{1/4} \right\} \ll \frac{x^{1/2 - \varepsilon}}{q^2}.
\end{equation*}
\end{corollary}

\subsection{Moment of order $1$ and correlations}

The Cauchy-Schwarz inequality immediately gives the following upper bound for the moment of order $1$ of $|E(x;q,a)|$.

\begin{corollary}
\label{Corollary 4}
Let $\varepsilon > 0$ be fixed. For $q \leq x$, we have
\begin{equation*}
\sideset{}{^\ast}\sum_{a \imod{q}} |E(x;q,a)| \ll x^{\varepsilon} \left( x^{1/4} q^{3/4} + x^{1/2} q^{1/4} \right).
\end{equation*}
\end{corollary}

In addition, for $\gamma$ a bijection of $(\mathbb{Z} / q \mathbb{Z})^{\times}$, let us define
\begin{equation*}
T_{\gamma}(x;q) =
\sideset{}{^\ast}\sum_{\substack{n_1, n_2 \leq x \\ n_1 = \gamma(n_2) \imod{q}}} |\mu(n_1)| |\mu(n_2)|.
\end{equation*}
Using the Cauchy-Schwarz inequality and Theorem \ref{Main Theorem}, we obtain an upper bound for the correlations of $E(x;q,a)$ and $E(x;q,\gamma(a))$. We will prove that this upper bound implies that $T_{\gamma}(x;q)$ satisfies the same asymptotic formula as $T(x;q)$.

\begin{corollary}
\label{Corollary 5}
Let $\varepsilon > 0$ be fixed. Let $\gamma$ be a bijection of $(\mathbb{Z} / q \mathbb{Z})^{\times}$. For $q \leq x$, we have
\begin{equation*}
T_{\gamma}(x;q) - \frac1{\varphi(q)} \left( \sideset{}{^\ast}\sum_{n \leq x} |\mu(n)| \right)^2 \ll 
x^{\varepsilon} \left( x^{1/2} q^{1/2} + \frac{x}{q^{1/2}} \right),
\end{equation*}
uniformly in $\gamma$.
\end{corollary}

\subsection{Acknowledgements}

It is a great pleasure for the author to thank Philippe Michel and Ramon Moreira Nunes for interesting conversations related to the topics of this article.

The financial support and the perfect working conditions provided by the \'{E}cole Polytechnique F\'{e}d\'{e}rale de Lausanne are gratefully acknowledged.

\section{Preliminaries}

The proof of Proposition \ref{Proposition} makes use of geometry of numbers results that we record here. The following lemma is due to Heath-Brown \cite[Lemma $3$]{MR757475} and provides an upper bound for the number of primitive integral solutions to a linear equation in three variables.

\begin{lemma}
\label{Lemma 1}
Let $(w_0,w_1,w_2) \in \mathbb{Z}^3$ be a primitive vector and let $U_0, U_1, U_2 \geq 1$. The number of primitive vectors $(u_0,u_1,u_2) \in \mathbb{Z}^3$ satisfying $|u_i| \leq U_i$ for $i \in \{ 0, 1, 2 \}$ and the equation
\begin{equation*}
u_0 w_0 + u_1 w_1 + u_2 w_2 = 0,
\end{equation*}
is at most
\begin{equation*}
12 \pi \frac{U_0 U_1 U_2}{\max \{ |w_i| U_i \} } + 4,
\end{equation*}
where the maximum is taken over $i \in \{ 0, 1, 2 \}$.
\end{lemma}

The next result is due to the author \cite[Lemma $2$]{MR3263143}. It gives an asymptotic formula for the number of solutions to certain linear congruences in two variables.

\begin{lemma}
\label{Lemma 2}
Let $\varepsilon > 0$ be fixed. Let $V_1, V_2 \geq 1$, and let $q \geq 1$ and $a_1, a_2 \in \mathbb{Z}_{\neq 0}$ be coprime to
$q$. Let also
\begin{equation*}
N(V_1,V_2;q,a_1,a_2) = \# \left\{ (v_1,v_2) \in \left( [1,V_1] \times [1,V_2] \right) \cap \mathbb{Z}^2,
\begin{array}{l}
a_1 v_1 = a_2 v_2 \imod{q} \\
\gcd(v_1v_2,q) = 1
\end{array}
\right\},
\end{equation*}
and
\begin{equation*}
N^{\ast}(V_1,V_2;q) = \frac1{\varphi(q)} \# \left\{ (v_1,v_2) \in \left( [1,V_1] \times [1,V_2] \right) \cap \mathbb{Z}^2, \gcd(v_1v_2,q) = 1 \right\}.
\end{equation*}
We have the estimate
\begin{equation*}
N(V_1,V_2;q,a_1,a_2) - N^{\ast}(V_1,V_2;q) \ll q^{\varepsilon} M(q,a_1,a_2),
\end{equation*}
where
\begin{equation*}
M(q,a_1,a_2) = \sum_{d \mid q} d \sum_{\substack{0 < |r|, |s| \leq q/2 \\ a_1 s + a_2 r = 0 \imod{d}}} |r|^{-1} |s|^{-1}.
\end{equation*}
\end{lemma}

The next lemma is also due to the author \cite[Lemma $9$]{MR3263143}. It provides an upper bound for the average of the quantity $M(q,f_1^2,f_2^2)$ while summed over coprime variables $f_1$ and $f_2$.

\begin{lemma}
\label{Lemma 3}
Let $\varepsilon > 0$ be fixed and $F_1, F_2 \geq 1/2$. We have the bound
\begin{equation*}
\sideset{}{^\ast}\sum_{\substack{F_i < f_i \leq 2F_i \\ \gcd(f_1,f_2)=1}} M(q,f_1^2,f_2^2) \ll
q^{\varepsilon} \left( F_1 F_2 + q \right),
\end{equation*}
where the symbol $\ast$ indicates that the summation is restricted to integers which are coprime to $q$, and where $i$ implicitly runs over the set $\{1,2\}$.
\end{lemma}

It is worth noting that Lemmas \ref{Lemma 2} and \ref{Lemma 3} are not as precise as the original results
\cite[Lemmas $2$ and $9$]{MR3263143} but the versions recorded here are sufficient for our purpose.

\section{Proofs}

\subsection{Proof of Proposition \ref{Proposition}}

\label{Section Proposition}

We have
\begin{equation*}
T(x;q) = \sideset{}{^\ast}\sum_{\substack{d_1 e_1^2, d_2 e_2^2 \leq x \\ d_1 e_1^2 = d_2 e_2^2 \imod{q}}} \mu(e_1) \mu(e_2).
\end{equation*}
We set $e = \gcd(e_1, e_2)$ and we write $e_i = e f_i$ for $i \in \{1,2\}$, where $f_1, f_2 \in \mathbb{Z}_{\geq 1}$ satisfy $\gcd(f_1,f_2) = 1$. If $d_1 f_1^2 \neq d_2 f_2^2$ then the fact that $e$ is coprime to $q$ implies that $e^2$ divides an integer $k \in \mathbb{Z}_{\neq 0}$ such that $|k| \leq xq^{-1}$, and thus $e \leq x^{1/2} q^{- 1/2}$. Moreover, let us bound the contribution of the $(d_1, d_2, e, f_1, f_2) \in \mathbb{Z}_{\geq 1}^5$ such that
$d_1 f_1^2 = d_2 f_2^2$ and $e > x^{1/2} q^{- 1/2}$. Since $\gcd(f_1,f_2) = 1$, there exists an integer
$m \in \mathbb{Z}_{\geq 1}$ such that $d_1 = m f_2^2$ and $d_2 = m f_1^2$. Therefore, this contribution is at most
\begin{equation*}
\sum_{\substack{m e^2 f_1^2 f_2^2 \leq x \\ e > x^{1/2} q^{-1/2}}} 1 \ll x^{1/2} q^{1/2}.
\end{equation*}
As a result, we have
\begin{equation}
\label{Error 1}
T(x;q) - T'(x;q) \ll x^{1/2} q^{1/2},
\end{equation}
where
\begin{equation*}
T'(x;q) = \sideset{}{^\ast}\sum_{e \leq x^{1/2} q^{-1/2}}
\ \sideset{}{^\ast}\sum_{\substack{f_1, f_2 \leq x^{1/2} e^{-1} \\ \gcd(f_1,f_2) = 1}}
\mu(e f_1) \mu(e f_2) N(D_1,D_2;q,f_1^2,f_2^2),
\end{equation*}
where we have set $D_i = xe^{-2}f_i^{-2}$ for $i \in \{1,2\}$, and where $N(D_1,D_2;q,f_1^2,f_2^2)$ is defined in Lemma \ref{Lemma 2}.

Let $y \geq 1$ be a parameter to be selected in due course. We define 
\begin{equation*}
T'_{\leq}(x;q) = \sideset{}{^\ast}\sum_{e \leq x^{1/2} q^{-1/2}}
\ \sideset{}{^\ast}\sum_{\substack{f_1, f_2 \leq x^{1/2} e^{-1} \\ \gcd(f_1,f_2) = 1 \\ ef_1f_2 \leq y}}
\mu(e f_1) \mu(e f_2) N(D_1,D_2;q,f_1^2,f_2^2).
\end{equation*}
We have
\begin{align*}
T'(x;q) - T'_{\leq}(x;q) & \ll \sideset{}{^\ast}\sum_{e^2 d_1 d_2 \leq x^2/y^2} \ 
\sideset{}{^\ast}\sum_{\substack{d_1 e^2 f_1^2, d_2 e^2 f_2^2 \leq x \\ d_1 f_1^2 = d_2 f_2^2 \imod{q} \\ \gcd(f_1,f_2) = 1}} 1
\\
& \ll \sideset{}{^\ast}\sum_{e^2 d_1 d_2 \leq x^2/y^2} \
\sum_{\substack{1 \leq \rho \leq q \\ \rho^2 d_1 = d_2 \imod{q}}} \ 
\sideset{}{^\ast}\sum_{\substack{d_1 e^2 f_1^2, d_2 e^2 f_2^2 \leq x \\ f_1 = \rho f_2 \imod{q} \\
\gcd(f_1,f_2) = 1}} 1.
\end{align*}
We now use Lemma \ref{Lemma 1}, we get
\begin{equation*}
\sideset{}{^\ast}\sum_{\substack{d_1 e^2 f_1^2, d_2 e^2 f_2^2 \leq x \\ f_1 = \rho f_2 \imod{q} \\
\gcd(f_1,f_2) = 1}} 1 \ll
\frac{x}{e^2 d_1^{1/2} d_2^{1/2} q} + 1.
\end{equation*}
Therefore, we find that
\begin{equation}
\label{Error 2}
T'(x;q) - T'_{\leq}(x;q) \ll x^{\varepsilon} \left( \frac{x^2}{y q} + \frac{x^2}{y^2} \right).
\end{equation}

We now estimate the quantity $T'_{\leq}(x;q)$ by making use of Lemma \ref{Lemma 2}. We obtain
\begin{equation}
\label{Error 3}
N(D_1,D_2;q,f_1^2,f_2^2) - N^{\ast}(D_1,D_2;q) \ll q^{\varepsilon} M(q,f_1^2,f_2^2),
\end{equation}
where $N^{\ast}(D_1,D_2;q)$ and $M(q,f_1^2,f_2^2)$ are defined in Lemma \ref{Lemma 2}. We now need to prove an upper bound for
\begin{equation*}
\sum_{e \leq x^{1/2} q^{-1/2}} \ \sideset{}{^\ast}\sum_{\substack{e f_1 f_2 \leq y \\ \gcd(f_1,f_2) = 1}}
M(q,f_1^2,f_2^2).
\end{equation*}
We split the summations over $f_1$ and $f_2$ into dyadic ranges. Let $F_1, F_2 \geq 1/2$ run over the set of powers of $2$.
Lemma \ref{Lemma 3} states that
\begin{equation*}
\sideset{}{^\ast}\sum_{\substack{F_i < f_i \leq 2 F_i \\ \gcd(f_1,f_2) = 1}} M(q,f_1^2,f_2^2) \ll
q^{\varepsilon} \left( F_1 F_2 + q \right).
\end{equation*}
We thus obtain that the overall contribution of the error term coming from the use of Lemma \ref{Lemma 3} is at most
\begin{equation*}
q^{\varepsilon} \sum_{e \leq x^{1/2} q^{-1/2}} \ \sum_{e F_1 F_2 \leq y} \left( F_1 F_2 + q \right) \ll
x^{\varepsilon} \left( y + x^{1/2} q^{1/2} \right).
\end{equation*}

As a result, recalling the upper bounds \eqref{Error 1}, \eqref{Error 2} and \eqref{Error 3}, we get
\begin{equation}
\label{Conclusion}
T(x;q) -
\frac1{\varphi(q)} \ \sideset{}{^\ast}\sum_{\substack{d_1 e^2 f_1^2, d_2 e^2 f_2^2 \leq x \\ \gcd(f_1,f_2) = 1 \\
e \leq x^{1/2} q^{-1/2} \\ e f_1 f_2 \leq y}} \mu(e f_1) \mu(e f_2) \ll
x^{\varepsilon} \left( x^{1/2} q^{1/2} + \frac{x^2}{y q} + \frac{x^2}{y^2} + y \right).
\end{equation}
Removing the conditions $e \leq x^{1/2} q^{-1/2}$ and $e f_1 f_2 \leq y$ from the summation in the estimate \eqref{Conclusion} creates an error term which is at most
$x^{\varepsilon} \left( x^{1/2} q^{1/2} + x^2 y^{-1} q^{-1} \right)$. In addition, we have
\begin{align*}
\sideset{}{^\ast}\sum_{\substack{d_1 e^2 f_1^2, d_2 e^2 f_2^2 \leq x \\ \gcd(f_1,f_2) = 1}} \mu(e f_1) \mu(e f_2) & = \sideset{}{^\ast}\sum_{d_1 e_1^2, d_2 e_2^2 \leq x} \mu(e_1) \mu(e_2) \\
& = \left( \sideset{}{^\ast}\sum_{n \leq x} |\mu(n)| \right)^2.
\end{align*}
We have finally obtained
\begin{equation*}
T(x;q) - \frac1{\varphi(q)} \left( \sideset{}{^\ast}\sum_{n \leq x} |\mu(n)| \right)^2 \ll
x^{\varepsilon} \left( x^{1/2} q^{1/2} + \frac{x^2}{y q} + \frac{x^2}{y^2} + y \right).
\end{equation*}
Choosing $y = x q^{-1/2}$ completes the proof of Proposition \ref{Proposition}.

\subsection{Proof of the equivalence between Theorem \ref{Main Theorem} and Proposition \ref{Proposition}}

We have
\begin{equation*}
V(x;q) = T(x;q) - 2 c_q \frac{x}{q} \sideset{}{^\ast}\sum_{n \leq x} |\mu(n)| + \varphi(q) c_q^2 \frac{x^2}{q^2},
\end{equation*}
and thus
\begin{equation*}
V(x;q) = T(x;q) - \frac1{\varphi(q)} \left( \sideset{}{^\ast}\sum_{n \leq x} |\mu(n)| \right)^2 +
\frac1{\varphi(q)} \left( \sideset{}{^\ast}\sum_{n \leq x} |\mu(n)| - \varphi(q) c_q \frac{x}{q} \right)^2.
\end{equation*}
Moreover, it is immediate to check that
\begin{equation*}
\sideset{}{^\ast}\sum_{n \leq x} |\mu(n)| - \varphi(q) c_q \frac{x}{q} \ll x^{1/2 + \varepsilon},
\end{equation*}
so Theorem \ref{Main Theorem} and Proposition \ref{Proposition} are clearly seen to be equivalent.

\subsection{Proof of Corollary \ref{Corollary 5}}

Let
\begin{equation*}
V_{\gamma}(x;q) = \sideset{}{^\ast}\sum_{a \imod{q}} E(x;q,a) E(x;q,\gamma(a)).
\end{equation*}
We have
\begin{equation*}
V_{\gamma}(x;q) = T_{\gamma}(x;q) - 2 c_q \frac{x}{q} \sideset{}{^\ast}\sum_{n \leq x} |\mu(n)| +
\varphi(q) c_q^2 \frac{x^2}{q^2},
\end{equation*}
and thus
\begin{equation*}
T(x;q) - T_{\gamma}(x;q) = V(x;q) - V_{\gamma}(x;q).
\end{equation*}
The Cauchy-Schwarz inequality gives
\begin{equation*}
|V_{\gamma}(x;q)| \leq V(x;q).
\end{equation*}
Therefore, we have $ T(x;q) - T_{\gamma}(x;q) \geq 0$ and
\begin{equation*}
T(x;q) - T_{\gamma}(x;q) \leq 2 V(x;q).
\end{equation*}
Corollary \ref{Corollary 5} is now a direct consequence of Theorem \ref{Main Theorem} and Proposition \ref{Proposition}.

\bibliographystyle{amsalpha}
\bibliography{biblio}

\end{document}